\documentclass{amsart}

\theoremstyle{definition}

\theoremstyle{remark}

\newtheorem{Definition}{\bf Definition}[section]
\newtheorem{Thm}[Definition]{\bf Theorem}
\newtheorem{Lem}[Definition]{\bf Lemma}
\newtheorem{Cor}[Definition]{\bf Corollary}

\newtheorem{Note}[Definition]{\bf Note}

\numberwithin{equation}{section}

\newcommand{\abs}[1]{\lvert#1\rvert}

\reversemarginpar

\newcommand{\nn}{\nonumber}
\newcommand{\no}{\noindent}
\newcommand{\realpart}{\mathop{\rm Re}\nolimits}

\newcommand{\ba}{\begin{eqnarray}}
\newcommand{\ea}{\end{eqnarray}}
\newcommand{\ift}{\int_{0}^{\infty}}
\newcommand{\ione}{\int_{0}^{1}}

\newcommand{\allR}{\mathbb{R}}
\newcommand{\allC}{\mathbb{C}}
\newcommand{\allN}{\mathbb{N}}
\newcommand{\nnN}{\mathbb{N}_{0}}

\begin{document}

\title[Generalized polygamma function] {A generalized polygamma function}

\author{Olivier Espinosa}
\address{Departamento de F\'{\i}sica,
Universidad T\'{e}c. Federico Santa Mar\'{\i}a, Valpara\'{\i}so, Chile}
\email{olivier.espinosa@fis.utfsm.cl}

\author{Victor H. Moll}
\address{Department of Mathematics,
Tulane University, New Orleans, LA 70118}
\email{vhm@math.tulane.edu}

\subjclass{Primary 33}

\date{\today}


\begin{abstract}
We study the properties of a function $\psi(z, q)$ (the
generalized polygamma function), intimately connected with the
Hurwitz zeta function and defined for complex values of the
variables $z$ and $q$, which is entire in the variable $z$ and
reduces to the usual polygamma function $\psi^{(m)}(q)$ for $z$ a
non-negative integer $m$, and to the balanced negapolygamma
function $\psi^{(-m)}(q)$ introduced in \cite{esmo2} for $z$ a
negative integer $-m$.
\end{abstract}

\maketitle


\section{Introduction} \label{S:intro}

The Hurwitz zeta function defined by
\ba
\zeta(z,q) & =  & \sum_{n=0}^{\infty} \frac{1}{(n+q)^{z}}
\label{zeta-def}
\ea
\no
for $z \in \allC, \; \realpart{z} > 1$ and $q \neq 0, \, -1, \, -2, \ldots$
is a generalization of the Riemann zeta function $\zeta(z) = \zeta(z,1)$.
This function admits a meromorphic continuation into the whole complex plane.
The only singularity is a simple pole at $z=1$ with unit residue.
The recent paper  \cite{knopp}
presents a motivated discussion of this extension.
\medskip

The Hurwitz zeta function turns out to be related to the classical
gamma function, defined for $\realpart{q} > 0$ by
\ba
\Gamma(q) & = & \ift t^{q-1} e^{-t} \, dt,
\label{gamma-def}
\ea
in several different ways. For example, the {\em digamma} function
\ba
\psi(q) & = & \frac{d}{dq} \log \Gamma(q)
\label{digamma-def}
\ea
\no
appears in the Laurent expansion of $\zeta(z,q)$ at the pole $z=1$:
\ba
\zeta(z,q) & = & \frac{1}{z-1} - \psi(q) + O(z-1).
\label{laurent}
\ea
\no
A second connection among these functions is given by Lerch's identity
\ba
\zeta'(0,q) & = & \log \Gamma(q) + \zeta'(0) = \log \left( \frac{\Gamma(q)}
{\sqrt{2 \pi}} \right)
\label{lerch}
\ea
\no
where we have used the classical value $\zeta'(0) = - \log \sqrt{2 \pi}$ in
the last step.  \\

A third example is the relation between the Hurwitz zeta function
and the {\em poly\-gam\-ma function} defined by
\ba
\psi^{(m)}(q) = \frac{d^{m}}{dq^{m}}
\psi(q), \quad m \in \allN,
\label{polygamma-def}
\ea
\no
namely
\ba
\psi^{(m)}(q) & = & (-1)^{m+1} m! \, \zeta(m+1,q).
\label{polygamma-hurwitz}
\ea
\no
These relations are not independent. Both (\ref{lerch}) and
(\ref{polygamma-hurwitz}) can be derived from (\ref{laurent}), in the 
limiting case $z \to 1$,  with the aid
of the formula
\ba
\left(\frac{\partial}{\partial q} \right)^{m} \zeta(z,q) & = & (-1)^{m} (z)_{m} 
 \, \zeta(z+m,q).
\label{hurwitz-der}
\ea

The digamma ($\psi(q) = \psi^{(0)}(q)$) and polygamma functions are analytic
everywhere in the complex $q$-plane, except for poles (of order
$m+1$) at all non-positive integers. The residues at these poles
are all given by $(-1)^{m+1} m!$.

\medskip

Extensions of the polygamma function $\psi^{(m)}(q)$ for $m$ a negative integer
have been defined by several authors \cite{adamchik, gosper, esmo2}. These
functions have been called {\em negapolygamma} functions. For example,
Gosper \cite{gosper} defined the family of functions
\begin{equation}\label{gosper-negapoly}
\begin{split}
\psi_{-1}(q) & := \log \Gamma(q), \\
\psi_{-k}(q) & := \int_0^q \psi_{-k+1}(t) dt,\quad k\ge 2,\\
\end{split}
\end{equation}
which were later reconsidered by Adamchik \cite{adamchik} in the form
\ba\label{adam-negapoly}
\psi_{-k}(q) = \frac{1}{(k-2)!}\int_0^q (q-t)^{k-2}\log\Gamma(t) dt,\quad k\ge
2.
\ea
These negapolygamma functions can be expressed in terms of
the derivative (with respect to its first argument) of the Hurwitz
zeta function at the negative integers
\cite{adamchik,gosper}.
\no
The definition of the negapolygamma functions in
(\ref{gosper-negapoly}) can be modified by introducing
arbitrary constants of integration at every step. This yields
infinitely many different families of negapolygamma functions,
with the property that the corresponding members of any two
families differ by a polynomial,
\[
\psi _a^{(-m)} (q) - \psi _b^{(-m)} (q) = p_{m - 1} (q),
\]
where the functions $p_n (q)$ are polynomials in $q$ of degree $n$,
satisfying the property
\[
p_n (q) = \frac{d}{{dq}}p_{n + 1} (q).
\]

An example of such modified negapolygamma functions has been
introduced in \cite{esmo2}, in connection with integrals involving
the polygamma and the loggamma functions.
These are the {\em balanced negapolygamma} functions, defined for $m \in
\allN$ by
\\
\ba\label{bal-negapolygamma}
\psi ^{( -m)} (q): = \frac{1}{m!}\left[{A_m (q)} - H_{m-1} B_m(q)\right],
\ea
\\
where $H_{r}:= 1 + 1/2 + \cdots + 1/r$ is
the harmonic number ($H_0:=0$), $B_m(q)$ is the $m$-th Bernoulli
polynomial, and the functions $A_m(q)$ are defined in terms of the
Hurwitz zeta function as
\ba\label{Ak-def}
A_{m}(q) :=  m\,\zeta'(1-m,q).
\ea
A function $f(q)$ is defined to be {\em balanced} (on the unit
interval) if it satisfies the two properties
\[
\ione f(q)dq=0 \quad \text{and}\quad f(0)=f(1).
\]
Note that the Bernoulli polynomials, which are related to the Hurwitz
zeta function in a way similar to (\ref{Ak-def}),
\ba\label{bernoulli-hurwitz}
B_{m}(q) = -m \, \zeta(1-m,q),  \quad m \in \allN,
\ea
are themselves balanced functions. In
\cite{esmo2} we have shown that the balanced negapolygamma functions
(\ref{bal-negapolygamma}) satisfy
\ba\label{der-negapolygamma}
\frac{d}{dq}\psi^{(-m)}(q)=\psi^{(-m+1)}(q),\quad m\in\allN.
\ea
This makes them a negapolygamma family, connecting
$\psi^{(-1)}(q)= \log\Gamma(q)+\zeta'(0)$ to the digamma function
$\psi^{(0)}(q)=d\log\Gamma(q)/dq$.
\\

The goal of this work is to introduce and study a meromorphic
function of two complex variables, $\psi(z,q)$, the {\em
generalized polygamma function}, that reduces to the polygamma
function $\psi^{(m)}(q)$ for $z = m \in \nnN$ and to the balanced
negapolygamma function $\psi^{(-m)}(q)$ for $z = - m \in - \allN$.
We describe some analytic properties of $\psi(z,q)$ and show they
extend those of polygamma and balanced negapolygamma functions. We
also present some definite integral formulas involving $\psi(z,q)$
in the integrand. Finally, we compare our generalized polygamma
function with a different generalization introduced by Grossman in
\cite{grossman}.\\

\medskip
\section{The generalized polygamma function}
\label{S:generalized-polygamma}

The generalized polygamma function is defined by
\ba
\label{psi-def-a}
\psi(z,q) := e^{- \gamma z } \frac{\partial}{\partial z}
\left[ e^{\gamma z} \frac{\zeta(z+1,q) }{\Gamma(-z)} \right],
\ea
where $z \in \allC$ and $q \in \allC, q\not\in -\nnN$.
At $z=m\in\allN$, where $\Gamma(-z)$ has a pole, and at $z=0$,
where both $\Gamma(-z)$ and $\zeta(z+1,q)$ have poles, we define
(\ref{psi-def-a}) by its corresponding limiting values given in the proof of 
Theorem \ref{psi-reduces-to-polygamma-thm}. We show 
below that, for fixed $q$, $\psi(z,q)$ is indeed an \emph{entire}
function of $z$.

The alternative representation
\ba
\psi(z,q) & = & e^{- \gamma z} \frac{\partial}{\partial z}
\left[ \frac{e^{\gamma z} }{\Gamma(1-z)}
\frac{\partial \zeta(z,q)}{\partial q} \right]
\label{pis-alter}
\ea
follows directly from (\ref{hurwitz-der}). \\

\begin{Lem}\label{psi-representations-lemma}
The function $\psi(z,q)$ is given by
\begin{align}
\label{psi-rep-b}
\psi(z,q)  & = \frac{1}{\Gamma(-z)}
\Big[
\zeta'(z+1,q)+\left\{\gamma+\psi(-z)\right\}\zeta(z+1,q)
\Big] \\
\intertext{and}
\label{psi-rep-c}
\psi(z,q)   & = \frac{1}{\Gamma(-z)}
\Big[
\zeta'(z+1,q)+H(-z-1)\zeta(z+1,q)
\Big],
\end{align}
where $H$ is defined by
\ba\label{harmonic-def}
H(z) := \sum_{k=1}^{\infty} \left( \frac{1}{k} - \frac{1}{k+z} \right).
\ea
\end{Lem}
\begin{proof}
Differentiation of (\ref{psi-def-a}) yields
(\ref{psi-rep-b}).  The second
representation follows from the identity $H(z) = \gamma +
\psi(z+1)$; see \cite{gr}, for instance.
\end{proof}

\no
The function $H$ can be termed the generalized harmonic number
function. It  has simple poles with residue $-1$ at all negative
integers, and reduces to the $n$-th harmonic number $H_k$ for
$z=k\in\nnN$. It satisfies the following reflection formula:
\ba\label{harmonic-reflection-formula}
H(-z)=H(z-1)+\pi\cot\pi z.
\ea

\medskip

We show first that, for $m \in \allN$, $\psi(-m,q)$ reduces to the
balanced negapolygamma function $\psi^{(-m)}(q)$ defined in
(\ref{bal-negapolygamma}).

\begin{Thm}\label{psi-reduces-to-negapolygamma-thm}
For $m \in \allN$, $\psi(-m,q) = \psi^{(-m)}(q)$.
\end{Thm}
\begin{proof}
Lemma \ref{psi-representations-lemma} gives
\ba
\psi(-m,q) &  = & \frac{1}{\Gamma(m)}\left[\zeta'(1-m,q) +
H_{m-1}\zeta(1-m,q)\right].
\ea
\no
The result now follows from \eqref{bal-negapolygamma}, \eqref{Ak-def},
and \eqref{bernoulli-hurwitz}.\\
\end{proof}
\medskip

We show next that the generalized polygamma function has no
singularities in the complex $z$ plane and that $\psi(0,q)$ is
actually the digamma function $\psi(q)$. \\

\begin{Thm}\label{psi-main}
For fixed $q \in \allC$, the function $\psi(z,q)$ is an entire
function of $z$. Moreover $\psi(0,q) = \psi(q)$.
\end{Thm}
\begin{proof}
In the representation (\ref{psi-rep-c}), the term
$1/\Gamma(z)$ is entire and $\zeta(z,q)$
has only a simple pole at $z=1$ and is analytic for $z \neq 1$.
Thus $z=0$ is the only possible singularity for $\psi(z,q)$. This
singularity is removable because for $z$ near $0$
\begin{align}
\frac{\zeta'(z+1,q)}{\Gamma(-z)}  & =   \left( - \frac{1}{z^{2}} + O(z) 
\right) \times \left( -z + \gamma z^{2} + O(z^{3}) \right) \nn \\
 & =  \frac{1}{z} - \gamma  + O(z) \nn \\
\intertext{and}
\frac{H(-z-1) \, \zeta(z+1,q) }{\Gamma(-z)}  & =    -\frac{1}{z} + \gamma
+ \psi(q) + O(z), \nn
\end{align}
\no
so that $\psi(z,q) = \psi(q) + O(z)$.
\end{proof}

\medskip

We finally show that, for $m \in \allN$, $\psi(m,q)$ reduces to the
polygamma function $\psi^{(m)}(q)$ defined in
(\ref{polygamma-def}). \\

\begin{Thm}\label{psi-reduces-to-polygamma-thm}
The function $\psi(z,q)$ satisfies
\begin{align}
\frac{\partial}{\partial q} \psi(z,q) & =   \psi(z+1,q)
\label{psi-der} \\
\intertext{and}
\psi(m,q) & =   \psi^{(m)}(q),  \quad m \in \allN.
\label{psi-reduces-to-polygamma-eqn}
\end{align}
\end{Thm}
\begin{proof}
Use (\ref{hurwitz-der}) to produce
\ba
\frac{\partial}{\partial q} \psi(z,q) =
-e^{- \gamma z} \frac{\partial}{\partial z}
\left[ e^{\gamma z} \frac{(z+1) \, \zeta(z+2,q)}{\Gamma(-z)} \right] \nn
\ea
\no
and then use $\Gamma(-z) = -(z+1) \Gamma(-z-1)$ to obtain (\ref{psi-der}). \\

The identity (\ref{psi-reduces-to-polygamma-eqn}) follows by induction
from Theorem
\ref{psi-main} and (\ref{psi-der}), but we provide an alternative proof. Set
$z = m + \epsilon$ and consider (\ref{psi-rep-c}) as $\epsilon \to 0$. The
expansions
\ba
\frac{1}{\Gamma(-m-\epsilon)} = (-1)^{m+1} m! \, \epsilon +
O(\epsilon^{2}) & \text{ and } &
H(-1-m-\epsilon) = \frac{1}{\epsilon} + H_{m} + O(\epsilon)
\nn
\ea
\no
are the only terms that produce a nonvanishing contribution in
\eqref{psi-rep-c} as $\epsilon \to 0$. We conclude that
$\psi(m,q) = (-1)^{m+1} m! \, \zeta(m+1,q)$ and
the result follows from (\ref{polygamma-hurwitz}).
\end{proof}

\section{Functional relations}

The generalized polygamma function $\psi(z,q)$, as a function of $q$,
 satisfies some
simple algebraic and analytic relations. These are derived from
those of $\Gamma(z)$ and $\zeta(z,q)$. \\

\begin{Thm}\label{psi-shift-thm}
The function $\psi(z,q)$ satisfies
\ba\label{psi-shift}
\psi (z,q + 1) = \psi (z,q) + \frac{{\ln q - H(- z-1)}}
{{q^{z+1} \Gamma ( - z)}}.
\ea
\end{Thm}
\begin{proof}
The identity
\ba\label{Hurwitz-shift}
\zeta (z,q + 1) = \zeta (z,q) - \frac{1}
{{q^z }},
\ea
produces
\[
\psi (z,q + 1) = \psi (z,q) - e^{ - \gamma z} \frac{\partial }
{{\partial z}}\left[ {e^{\gamma z} \frac{1}
{{q^{z+1} \Gamma ( - z)}}} \right].
\]
The result now follows from
$\gamma  + \psi(-z)=H(-z-1)$.
\end{proof}

\medskip

Relation (\ref{psi-shift}) generalizes the well-known
functional relations for the digamma and polygamma functions,
\ba
\psi (q + 1) & = & \psi (q) + \frac{1}{q},
\label{digamma-shift} \\
\psi ^{(m)} (q + 1) & = & \psi ^{(m)} (q) + \frac{{( - 1)^m m!}}
{{q^{m + 1} }},
\label{polygamma-shift}
\ea
\no
and the corresponding relation
\ba
\psi^{(-m)}(q+1) & = & \psi^{(-m)}(q)+\frac{q^{m-1}}{(m-1)!}[\ln
q-H_{m-1}]
\label{negapolygamma-shift}
\ea
\no
for the balanced negapolygamma function \cite{esmo2}.

\medskip

\no{\bf Note.}
We have been unable to find a generalization of the other well-known
functional relation for the polygamma function,
\ba\label{polygamma-reflection}
( - 1)^m \psi ^{(m)} (1 - q) = \psi ^{(m)} (q) + \frac{{d^m }}
{{dq^m }}\pi \cot \pi q.
\ea

\medskip

The next result establishes a multiplication formula for
$\psi(z,q)$. It generalizes the analogous result for the digamma
function, \cite{gr} (8.365.6). \\

\begin{Thm}\label{psi-multiplication-thm}
Let $k \in \allN$. Then,
\ba
\label{psi-multiplication-eqn}
k^{z + 1} \psi (z,k q) & = & \sum\limits_{j = 0}^{k - 1}
{\psi \left( {z,q + j/k} \right)}  - k^{z + 1} \ln k \, \frac{{\zeta (z+1,k q)}}
{{\Gamma ( - z)}} \\
& = &
\sum\limits_{j = 0}^{k - 1}
{\left[ {\psi \left( {z,q + j/k} \right) - \frac{{\ln k}}
{{\Gamma ( - z)}}\zeta \left( {z+1,q + j/k} \right)} \right]}.
\nn
\ea
\end{Thm}
\begin{proof}
Use the multiplication rule
\ba\label{hurwitz-multiplication}
k^z \zeta (z,kq) =
\sum\limits_{j = 0}^{k - 1} {\zeta \left( {z,q + j/k} \right)}
\ea
for the Hurwitz zeta function in the definition
(\ref{psi-def-a}) of $\psi(z,q)$.
\end{proof}
\no
The case $k=2$ yields the duplication formula
\ba\label{psi-duplication}
\psi (z,2q) = \frac{1}
{{2^{z + 1} }}\left[ {\psi (z,q) + \psi (z,q + 1/2)} \right] - \ln 2 \, \frac{{\zeta
(z+1,2q)}}
{{\Gamma ( - z)}}.
\ea

\section{Series expansions of $\psi (z,q)$}

In this section we present two different series expansions for the
generalized polygamma function. The first is a generalization of
the well-known expansion of the digamma function,
\ba\label{digamma-expansion}
\psi(q+1) = - \gamma + \sum_{k=1}^{\infty} (-1)^{k+1} \zeta(k+1) q^{k},
\quad |q| < 1.
\ea

\medskip

\begin{Thm}
\label{psi-series-thm}
Let $z \in \mathbb{C}$ and $|q| < 1$. Then
\ba\label{psi-series}
\psi(z,q+1) = \sum_{k=0}^{\infty} \psi(z+k,1) \, \frac{q^{k}}{k!}.
\ea
\end{Thm}
\begin{proof}
The Taylor expansion of $\psi(z,q+1)$ around $q=0$ can be
expressed in terms of
$\psi(z,q)$ using the iterated version of (\ref{psi-der}),
\ba\label{psi-higher-der}
\frac{\partial^{k}}{\partial q^{k}} \psi(z,q) = \psi(z+k,q)
\ea
\no
evaluated at $q=1$. The radius of convergence is computed to be
$1$ by using the ratio test, the identity
\ba\label{psi-at-one-1}
\psi(z,1) =  \frac{1}{\Gamma(-z)} \left[ \zeta'(z+1) + H(-z-1) \, \zeta(z+1)
\right],
\label{psiat1}
\ea
\no
and the fact that  $\zeta'(z+1)$ tends to zero faster than the
term $H(-z-1) \, \zeta(z+1)$ as $z\to\infty$.
\end{proof}


\medskip

\no
\begin{Note}\label{note-loggamma-expansion}
The series expansion (\ref{digamma-expansion}) for the digamma
function is the special case $z = 0$ of (\ref{psi-series}). This
follows from the values
\begin{align}
\label{psi-at-0-1}
\psi(0,1)&=\psi(1)=-\gamma
\intertext{and}
\label{psi-at-k-1}
\psi(k,1)&=\psi^{(k)}(1)=(-1)^{k + 1} k!\,\zeta (k + 1),
\end{align}
for $k \in \allN$. Similarly $z = -1$ and the value
$\psi(-1,1)=\zeta'(0)$ yield
the well-known expansion of the loggamma function,
\ba
\label{loggamma-expansion}
\log\Gamma(q+1) & =  &
- \gamma q + \sum_{k=2}^{\infty} (-1)^{k} \frac{\zeta(k)}{k} q^{k},
\quad |q|  < 1.
\ea
\end{Note}

\no
\begin{Note}
Riemann's functional equation,
\ba
\zeta(1-z)  =  \frac{\zeta(z) \, (2 \pi)^{1-z}}
{2 \Gamma(1-z) \, \sin(\pi z/2) }
\label{zeta-functional-1}
\\
= 2 \cos \left( \frac{\pi z}{2} \right) \, \frac{\zeta(z) \Gamma(z) }
{(2 \pi)^{z}}, \nn
\ea
\no
yields the alternate representation
\begin{multline}
\label{psi-at-one-2}
\quad
\psi(z,1) = 2 \, (2 \pi)^{z} \, \cos \left( \frac{\pi z}{2} \right)
\,\left[ \left( \gamma + \ln 2 \pi - \frac{\pi}{2} \tan \frac{\pi z}{2}
\right) \zeta(-z) - \zeta'(-z) \right].
\quad
\end{multline}
\end{Note}

\medskip

\begin{Note}
Theorems \ref{psi-shift-thm} and
\ref{psi-series-thm} determine the behavior of $\psi(z,q)$ for small $q$:
\begin{multline}\label{psi-at-zero}
\psi (z,q) =  - \frac{1}
{{\Gamma ( - z)}}\frac{{\ln q}}
{{q^{z+1} }} + \frac{{H(- z-1)}}
{{\Gamma ( - z)}}\frac{1}
{{q^{z+1} }} + \psi (z,1) + \psi (z + 1,1)q + O(q^2 ).
\end{multline}
\no
For $ z = m \in \nnN$ the coefficients of the first two terms are
\[
\frac{1}{\Gamma(-m)} = 0 \quad\text{ and }\quad
\frac{H(-m-1)}{\Gamma(-m)} = (-1)^{m+1} \, m!,
\]
\no
so the logarithmic term drops out and we recover the known
behavior of the polygamma function as $q \to 0$,
\ba
\psi^{(m)}(q) = \frac{(-1)^{m+1} \, m!}{q^{m+1}} + \psi^{(m)}(1) +
O(q).
\ea
\no
For $ z \not \in \allN$ with $\realpart{z} \geq -1$ the first term
in (\ref{psi-at-zero}) determines the leading behavior, and if
$\realpart{z} <  -1$ the first two terms in (\ref{psi-at-zero})
vanish  as $q \to 0$ and hence $\psi(z,q)$ tends to the finite
value $\psi(z,1)$ given by (\ref{psi-at-one-1}) or
(\ref{psi-at-one-2}).
\end{Note}

\medskip

\no
We now establish a Fourier series representation for the
generalized polygamma function $\psi(z,q)$.
\medskip

\begin{Thm}
For $\realpart z <  -1$ and $0\le q \le 1$:
\begin{multline}\label{psi-fourier}
\psi (z,q) = 2(2\pi )^z \left[ \sum\limits_{n = 1}^\infty
n^z (\gamma  + \ln 2\pi n)\cos (2\pi nq + \pi z/2)\right.\\
\left. - \frac{\pi }{2}\sum\limits_{n = 1}^\infty  n^z \sin (2\pi nq + \pi z/2)
\right].
\end{multline}
\end{Thm}
\no
This result generalizes the Fourier expansion for the
balanced nega\-poly\-gamma function given in \cite{esmo2}. It
implies that $\psi(z,q)$ is itself balanced for any $z$ such that
$\realpart{z} <  -1$.
\begin{proof}
Let $s = z+1$ in the Fourier representation of the Hurwitz zeta
function,
\begin{multline}
\label{hurwitz-fourier}
\zeta(s,q) = \frac{2 \Gamma(1-s)}{(2 \pi)^{1-s}}\;
\left[ \sin \left( \frac{\pi s}{2} \right) \sum_{n=1}^{\infty}
\frac{\cos( 2 \pi q n)}{n^{1-s}} +
\cos \left( \frac{\pi s}{2} \right)
\sum_{n=1}^{\infty} \frac{ \sin( 2 \pi q n) }{ n^{1-s} }
\right],
\end{multline}
which is valid for $\realpart s  <  0$ and $0 \le q \le 1$,
and substitute (\ref{hurwitz-fourier}) into
(\ref{psi-def-a}).
\end{proof}

\section{Integral representations of $\psi (z,q)$}
\label{S:integral-representations}

This section contains integral representations for $\psi(z,q)$
that are derived directly from corresponding integral
representations of the Hurwitz zeta function. For instance,
\ba\label{hurwitz-int-rep-1}
\zeta (z,q) = \frac{1}
{{\Gamma (z)}}\int_0^\infty  {\frac{{e^{ - qt} }}
{{1 - e^{ - t} }}\,t^{z-1} dt},
\ea
valid for $\realpart z > 1$ and $\realpart q > 0$, implies the next
result.

\medskip

\begin{Thm}
Let $\realpart{z}  > 0$ and $\realpart{q}  > 0$. Then
\ba
\psi(z,q) = - \ift \frac{e^{-qt} t^{z} }{1 - e^{-t}}
\left[ \cos \pi z + \frac{\gamma}{\pi} \sin \pi z + \frac{\sin \pi z}{\pi}
\, \ln t \right] \, dt.
\ea
\end{Thm}
\begin{proof}
The identity
\ba
\frac{\zeta(z+1,q)}{\Gamma(-z)} = - \frac{\sin \pi z}{\pi}
\ift \frac{e^{-qt}}{1 - e^{-t}} \, t^{z} \, dt
\ea
\no
follows from the integral representation for $\zeta(z,q)$ in
(\ref{hurwitz-int-rep-1}) and the reflection rule for the gamma
function. The result now follows from the definition of
$\psi(z,q)$.
\end{proof}

\medskip

A {\em Hankel-type} contour is a curve that starts at $+ \infty +
i \, 0+$, moves to the left on the upper half-plane, encircles the
origin once in the positive direction, and returns to $+\infty - i
\, 0+$ on the lower half-plane.
The Hurwitz zeta function has the following
integral representation along a Hankel-type
contour\cite{ww}:
\ba\label{hurwitz-int-rep-2}
\frac{{\zeta (z+1,q)}}
{{\Gamma ( - z)}} =  - \frac{1}
{{2\pi i}}\int_\infty ^{(0 + )} {\frac{{ e^{ - qt} }}
{{1 - e^{ - t} }}( - t)^z\,dt},
\ea
valid for arbitrary complex $z$ and $\realpart q  >  0$. \\
\no

\begin{Thm}
Let $q, \, z \in \allC$ with $\realpart{q} > 0$. Then
\ba\label{psi-int-rep-2}
\psi(z,q) = - \frac{1}{2 \pi i} \int_{\infty}^{(0+)}
\frac{[ \gamma + \ln(-t) ] e^{-qt} }{1 - e^{-t}} ( - t)^z\, dt.
\ea
\end{Thm}
\begin{proof}
The result follows directly from (\ref{hurwitz-int-rep-2}).
\end{proof}

\section{Definite integrals involving $\psi (z,q)$}
\label{S:integrals}

Definite integrals of $\psi (z,a+bq)$ can be directly obtained from
its primitive,
\ba\label{psi-int-0}
\int {\psi (z,a+bq)}\,dq
= b^{-1}\psi (z-1,a+bq),
\ea
according to (\ref{psi-der}). So, for example,
\ba\label{psi-int-1}
\int_0^1 {\psi (z,q)}\,dq
=\begin{cases}
0,  & \text{ if } \realpart{z}  <  0,  \\
\infty, & \text{ if } \realpart{z}  \ge 0,
\end{cases}
\ea
where we have used the result (\ref{psi-at-zero}) to evaluate
$\psi(z,q)$ at the origin.

\medskip

The integral formulas presented below are direct consequences of
the corresponding integral formulas for the Hurwitz
zeta function. Several of these were derived in
\cite{esmo1, esmo2}.
\begin{Thm}\label{psi-int-2-thm}
For $\realpart z, \realpart z' < 0$ and $\realpart (z+z') <  -1$,
\begin{multline}
\label{psi-int-2}
\int_0^1 {\psi (z,q)\psi (z',q)}\,dq
=
2(2\pi )^{z + z'} \cos \frac{{\pi (z - z')}}
{2}\Bigg\{ \left[ \frac{{\pi ^2 }}{4} +
\left( {\gamma  + \ln 2\pi } \right)^2 \right]
\zeta ( - z - z')\\
- 2\left( {\gamma  + \ln 2\pi } \right)\zeta '( - z - z')
+ \zeta ''( - z - z') \Bigg\}.
\end{multline}
\end{Thm}
\begin{proof}
This is a direct consequence of the following result \cite{esmo1},
\ba
\nn
\ione \zeta(s,q) \zeta(s',q) dq  =\frac{2 \Gamma(1-s) \Gamma(1-s')}
{(2 \pi)^{2 - s - s'}} \zeta(2 - s - s')
\cos \frac{\pi(s-s')}{2},
\ea
valid for $\realpart s < 1, \realpart s' < 1, \realpart(s+s') < 1$. Set
$s=z+1, s'=z'+1$, divide by $\Gamma(-z)\Gamma(-z')$, and construct
the functions $\psi(z,q), \psi(z',q)$ in the integrand according
to definition (\ref{psi-def-a}).
\end{proof}

\no
The evaluation (\ref{psi-int-2}) generalizes Example 5.6 in
\cite{esmo2}:
for $k,k'\in\allN$,
\begin{multline}\nn
\ione\psi ^{( - k)} (q)\psi ^{( - k')} (q)\,dq =
\frac{{2\cos (k - k')\frac{\pi}{2}}}{{(2\pi )^{k + k'} }}
\Bigg[ \zeta'' (k + k') - 2(\gamma  +
\ln 2\pi )\zeta' (k + k')\\
 + \left\{ {(\gamma  + \ln 2\pi )^2  +
\frac{{\pi ^2 }}{4}} \right\}\zeta (k + k') \Bigg].
\end{multline}
\no
The special case $k = k' =1$ reduces to
\ba
\int_{0}^{1} \left( \ln \Gamma(q) \right)^{2} dq &= &
\frac{\gamma^{2}}{12} + \frac{\pi^{2}}{48} + \frac{1}{3} \gamma
\ln \sqrt{2 \pi} + \frac{4}{3} \ln^{2} \sqrt{2 \pi} \nn \\
& &\phantom{xxx}-
(\gamma + 2 \ln \sqrt{2 \pi} ) \frac{\zeta'(2)}{\pi^{2}} +
\frac{\zeta''(2) }{2 \pi^{2}}, \nn
\ea
\no
given in \cite{esmo1}. \\

\begin{Cor} Let $\realpart z <  -1/2$. Then
\begin{align}
\label{psi-int-2-a}
&\int_0^1 {\psi (z,q)^2}\,dq
=
2(2\pi )^{2z}
\Bigg\{ \left[ \frac{{\pi ^2 }}{4} +
\left( {\gamma  + \ln 2\pi } \right)^2 \right]
\zeta ( - 2z)
\\
\nn
&\qquad\qquad\qquad\qquad\qquad\qquad
- 2\left( {\gamma  + \ln 2\pi } \right)\zeta '( - 2z)
+ \zeta ''( - 2z) \Bigg\}.
\end{align}
For $\realpart z <  -1$,
\ba
\label{psi-int-2-b}
\int_0^1 {\psi (z,q)\psi (z+1,q)}\,dq=0.
\ea
\end{Cor}

\medskip

\no
The same type of argument gives the next evaluation.  \\

\begin{Thm}
For $\realpart z, \realpart z' <  0$, and $\realpart (z+z') <  -1$,
\begin{multline}
\label{psi-int-3}
\int_0^1 \zeta (z+1,q)\psi (z',q)dq = 2(2\pi )^{z + z'}\Gamma (-z)
\Bigg\{
\frac{\pi }{2}\zeta (- z - z')\sin \frac{\pi }{2}(z - z')\\
+
\Big[ (\gamma  + \ln 2\pi )\zeta (-z - z')
- \zeta '(- z - z') \Big]\cos \frac{\pi }{2}(z - z')
\Bigg\}.
\end{multline}
\end{Thm}

\begin{Cor}
For $\realpart z <  0$,
\ba
\label{psi-int-3-a}
\int_0^1 \zeta (z,q)\psi (z,q)dq=
-\frac{1}{2}(2\pi )^{2z}\Gamma (1-z)
\zeta (1 - 2z).
\ea
\end{Cor}

\medskip

Our final evaluation computes the Mellin transform of the
generalized polygamma function.

\begin{Thm}\label{psi-mellin}
Let $a,b\in\allR^+$, $\alpha,z\in\allC$ such that
$0 <  \realpart\alpha < \realpart z$. Then
\begin{multline}
\label{psi-int-4}
\int_0^\infty {q^{\alpha-1}\psi (z,a+bq)}\,dq
=
\frac{{b^{ - \alpha } \Gamma (\alpha )}} {{\sin \pi (z - \alpha
)}}\Big[ (\sin \pi z)\psi (z - \alpha ,a) \\
+ (\sin \pi \alpha
)\Gamma (z+1 - \alpha )\zeta (z+1 - \alpha ,a) \Big].
\end{multline}
\end{Thm}
\begin{proof}
Start from formula (2.3.1.1) of \cite{prudnikov},
\ba
\ift q^{\alpha -1} \zeta(s,a + bq) dq = b^{-\alpha} B( \alpha, s - \alpha)
\zeta(s - \alpha, a),
\nn
\ea
\no
valid for $a, b \in \allR^{+}$ and $ 0 <  \realpart(\alpha) < 
\realpart(s) -1$, set $s=z+1$, and use the definition
(\ref{psi-def-a}) of $\psi(z,q)$ to evaluate the
integral as
\ba
\nn
 - \frac{{b^{ - \alpha } \Gamma (\alpha )}}
{\pi }e^{ - \gamma z} \frac{\partial }
{{\partial z}}\left[ {e^{\gamma z} (\sin \pi z)
\Gamma (z+1 - \alpha )\zeta (z+1 - \alpha ,a)} \right].
\ea
The desired evaluation now follows from the
reflection formulas for the gamma and digamma functions,
\ba
\nn
\Gamma(1-x)\Gamma(x)=\frac{\pi}{\sin\pi x}
\quad\text{and}\quad
\psi(1-x)=\psi(x)+\pi\cot\pi x,
\ea
respectively.
\end{proof}
\no

\no
\begin{Note}
The special case $z = m\in\allN$ in Theorem \ref{psi-mellin}
yields an explicit form for the Mellin transform of the
polygamma function $\psi^{(m)}(a + bq)$:
\begin{multline}
\label{polygamma-mellin}
\int_0^\infty {q^{\alpha-1}\psi^{(m)}(a+bq)}\,dq
=
( - 1)^{m + 1} b^{ - \alpha } \Gamma (\alpha )
\Gamma (1 + m - \alpha )\zeta (1 + m - \alpha,a),
\end{multline}
valid when $0 < \realpart\alpha < m$ and $a, \, b \in \mathbb{R}^{+}$.
This formula generalizes formula (6.473) of \cite{gr} to the
case $a,b\ne 1$.
\end{Note}

\section{Relation to Grossman's generalization of the polygamma
function}\label{S:other-generalizations}

In 1975, N. Grossman presented a generalization of
polygamma functions to arbitrary complex orders \cite{grossman}
which is different to ours.
He was motivated by a problem proposed a year earlier by B. Ross
\cite{ross} concerning the convergence and evaluation of the
integral
\ba\label{ross-integral}
I = \int_0^q (q-t)^{p-1}\log\Gamma(t) dt.
\ea
For $p\in\allN$, this integral corresponds precisely (up to a
normalization factor) to the Gosper-Adamchik's negapolygamma
functions defined by (\ref{adam-negapoly}). In \cite{grossman} the
author used the techniques of Liouville's fractional integration
and differentiation to obtain a generalization $\psi^{(\nu)}(q)$
of the polygamma function, with $\nu\in\allC$, in the form
\begin{multline}
\psi^{(\nu)}(q) = \frac{{q^{-\nu-1} }}
{{\Gamma (-\nu)}}\left\{ { \ln\frac{1}{q} + \gamma  + \frac{{\Gamma'(-\nu)}}
{{\Gamma (-\nu)}}} \right\}  - \frac{{\gamma q^{-\nu}}}
{{\Gamma (1-\nu)}}\\
- \frac{{q^{-\nu-1} }}
{{2\pi i}}\int_{\lambda  - i\infty }^{\lambda  + i\infty } {q^z \frac{{\Gamma
(z)\zeta (z)}}
{{\Gamma (z-\nu)}}\frac{\pi }
{{\sin \pi z}}\,} dz,
\end{multline}
where the contour of integration is along a vertical line with
$1 <  \lambda <  2$. The function $\psi^{(\nu)}(q)$ is an entire
function in the $\nu$-plane, for each $q$ in the plane cut along
the negative real axis \cite{grossman}. \\

For $\nu=-m\in -\nnN$, Grossman's generalized polygamma
$\psi^{(\nu)}(q)$ reduces to the Gosper-Adamchik
negapolygamma functions $\psi_{-m}(q)$. We showed in
\cite{esmo2} that the latter are related to the balanced
negapolygammas $\psi^{(-m)}(q)$ by
\begin{align}
\psi^{(-m)}(q) &=\psi_{-m}(q)+
\sum_{r=0}^{m-1}\frac{q^{m-r-1}}{r!(m-r-1)!}
\left[\zeta'(-r)+H_r\zeta(-r)\right],\\
\intertext{which, in light of (\ref{psi-at-one-1}), can be also
expressed as}
\label{negapolygammas-relationship}
\psi^{(-m)}(q) &=\psi_{-m}(q)+
\sum_{r=0}^{m-1}\frac{q^{m-r-1}}{\Gamma(m-r)}
\psi(-r-1,1).
\end{align}
\\
In the remainder of this section we shall explore the relation
between the functions $\psi(\nu,q)$ and $\psi^{(\nu)}(q)$ for
arbitrary values of the complex variable $\nu$. Since both of
these functions are entire in $\nu$, their difference
\ba\label{Psi-def}
\Psi(\nu,q):=\psi(\nu,q)-\psi^{(\nu)}(q)
\ea
must be an entire function itself. Furthermore, since both
$\psi(\nu,q)$ and $\psi^{(\nu)}(q)$ reduce to the standard
polygamma function when $\nu\in\nnN$, $\Psi(\nu,q)$ vanishes
identically at $\nu\in\nnN$. In order to study further properties
of the function $\Psi(\nu,q)$, we shall consider the asymptotic
and small-$q$ series expansions of both $\psi(\nu,q)$ and
$\psi^{(\nu)}(q)$.
First, we shall derive the correct asymptotic expansion of
Grossman's polygamma for large $q$, since this was incorrectly
given in \cite{grossman}. Let
\ba
I(\nu,q)= \frac{1}
{{2\pi i}}\int_{\lambda  - i\infty }^{\lambda  + i\infty } {q^z \frac{{\Gamma
(z)\zeta (z)}}
{{\Gamma (z-\nu)}}\frac{\pi }
{{\sin \pi z}}\,} dz.
\ea
As suggested in \cite{grossman},
for $ q > 1$ we can deform the contour so that it starts at
$-\infty-i0+$, runs below the real axis, encircles the point $z=1$
in the positive sense (crossing the real axis to the left of
$z=2$), and then returns to $-\infty+i0+$ running over the real
axis. $I(\nu,q)$ can then be evaluated along the deformed contour by a
residue calculation. The  only relevant poles are $z=1$ and $z=0,
-1, -2,\ldots$, coming from $\Gamma(z)$, $\zeta(z)$, and from the
zeros of $\sin(\pi z)$. The poles at $z = -2k$ are simple since
$\zeta(-2k)=0$. All the other poles are double.  Let
$R_\nu(z_0)$ be the residue at the pole $z=z_0$. Then
\begin{align}
R_\nu (1) & =  \frac{{ - q\ln q + q\psi (1 - \nu )}}
{{\Gamma (1 - \nu )}},
\nn\\
R_\nu  (0) & = \frac{{H( - 1 - \nu ) - \ln 2\pi  - \ln q}}
{{2\Gamma ( - \nu )}},
\nn\\
R_\nu  ( - m) & = \frac{1}
{{m!\Gamma ( - m - \nu )q^m }}\left[ {\zeta '( - m) - \frac{{B_{m + 1} }}
{{m + 1}}\left\{ {\ln q + H_m  - H( - m - \nu  - 1)} \right\}}
\right],
\nn
\end{align}
for $m=1,2,3,\ldots$. Using the
special values $B_0=1, B_1=-1/2, \zeta'(0)=-\frac{1}{2}\ln2\pi$, and
$H_0=0$, we obtain the asymptotic expansion
\begin{multline}
\psi ^{(\nu )} (q) \sim q^{ - \nu } \left\{ {\ln q\sum\limits_{k = 0}^\infty
{\frac{{B_k }}
{{k!\Gamma (1 - \nu  - k)q^k }}}
- \sum\limits_{k = 1}^\infty  {\frac{{k\zeta '(1 - k) - B_k H_{k - 1} }}
{{k!\Gamma (1 - \nu  - k)q^k }}} } \right.
\\
\left. {  - \sum\limits_{k = 0}^\infty  {\frac{{B_k H( - k -
\nu )}}
{{k!\Gamma (1 - \nu  - k)q^k }}} } \right\}.
\end{multline}
We observe that  Grossman \cite{grossman} missed most of the
logarithmic contribution. \\

The asymptotic expansion of the generalized
polygamma function $\psi(\nu,q)$ for large $q$ can be
obtained from \eqref{psi-def-a} and the asymptotic expansion of $\zeta(z,q)$
itself. This yields
\begin{multline}
\psi (\nu ,q) \sim q^{ - \nu } \left\{\ln q\frac{{\sin \pi \nu }}
{\pi }\sum\limits_{k = 0}^\infty  {( - 1)^k \frac{{B_k }}
{{k!}}\frac{{\Gamma (k + \nu )}}
{{q^{k } }}} - \cos \pi \nu \sum\limits_{k = 0}^\infty  ( - 1)^k \frac{{B_k }}
{{k!}}\frac{{\Gamma (k + \nu )}}{q^{k }}\right.
\\
\left. - \frac{{\sin \pi \nu }}
{\pi }\sum\limits_{k = 0}^\infty  {( - 1)^k \frac{{B_k }}
{{k!}}\frac{{H(k + \nu  - 1)\Gamma (k + \nu )}}
{{q^{k} }}}\right\}.
\end{multline}
The reflection formula for $\Gamma(z)$ yields
\[
\frac{{\sin \pi \nu }}
{\pi }( - 1)^k \Gamma (k + \nu ) = \frac{1}
{{\Gamma (1 - \nu  - k)}},
\]
and the reflection formula \eqref{harmonic-reflection-formula}
for the harmonic number function produces
\[
H(k+\nu-1)=H(-k-\nu)-\pi\cot\pi\nu.
\]
Thus
\ba
\psi (\nu ,q) \sim \ln q\sum\limits_{k = 0}^\infty  {\frac{{B_k }}
{{k!\Gamma (1 - \nu  - k)q^{k + \nu } }}}
-\sum\limits_{k = 0}^\infty  {\frac{{B_k H( - \nu  - k)}}
{{k!\Gamma (1 - \nu  - k)q^{k + \nu } }}}.
\ea
\no
We obtain therefore the following asymptotic expansion for the
function $\Psi(\nu, q)$ defined by \eqref{Psi-def}:
\ba
\Psi (\nu ,q) \sim \sum\limits_{k = 1}^\infty  {\frac{{\psi ( -
k,1)}}
{{\Gamma (1 - \nu  - k)q^{k + \nu } }}}.
\ea
We note that for $\nu=m\in\nnN$ the formal series on the
right-hand side vanishes identically, as it should.
For $\nu=-m\in-\allN$, the series above reduces to a
polynomial in $q$, which coincides with the one appearing in
\eqref{negapolygammas-relationship}. \\

On the other hand, for $\abs{q} < 1$, Grossman has proven that his
polygamma function has the convergent expansion\footnote{There is
actually an error in the expression given in \cite{grossman}.}
\begin{multline}
\psi ^{(\nu)} (q) = \frac{{q^{ - \nu  - 1} }}
{{\Gamma ( - \nu )}}\left\{ { - \ln q + \gamma  + \psi ( - \nu ) +
\frac{{\gamma q}}
{\nu } + \sum\limits_{k = 2}^\infty  {( - 1)^k \zeta (k)B( - \nu ,k)q^k } }
\right\},
\end{multline}
which, on account of the special values for $\psi(z,1)$ at the
non-negative integers given in \eqref{psi-at-0-1} and
\eqref{psi-at-k-1}, can be  written as
\ba
\psi^{(\nu)} (q) = \frac{{ - \ln q + H( - \nu  - 1)}}
{{q^{\nu  + 1} \Gamma ( - \nu )}} + \sum\limits_{k = 0}^\infty  {\frac{{\psi (k,1)}}
{{\Gamma ( - \nu  + k + 1)}}q^{k - \nu } }.
\ea
This is to be compared with the small-$q$ expansion of the
generalized polygamma function $\psi(\nu,q)$ obtained in Theorems
\ref{psi-shift-thm} and \ref{psi-series-thm}:
\ba
\psi (\nu ,q) = \frac{{ - \ln q + H( - \nu  - 1)}}
{{q^{\nu  + 1} \Gamma ( - \nu )}} + \sum\limits_{k = 0}^\infty  {\frac{{\psi (k +
\nu ,1)}}
{{\Gamma (k + 1)}}q^k }.
\ea
Again, both expansions coincide if $\nu\in\nnN$, since $1/\Gamma(z)$
vanishes at the non-positive integers, and differ by the
polynomial in \eqref{negapolygammas-relationship} if
$\nu=-m\in-\allN$.

\bigskip

\no
{\bf Acknowledgments}. The first author would like to thank the
Department of Mathematics at Tulane University for its hospitality
and acknowledge the partial support of grant MECESUP (Chile)
FSM-99-01. The second author acknowledges the partial support of
NSF-DMS 0070567.


\begin{thebibliography}{99}

\bibitem{adamchik}
ADAMCHIK, V.: {\em Polygamma functions of negative order}.
Jour. Comp. Appl. Math. {\bf 100}, 1998, 191-199.

\bibitem{boesmo1}
BOROS, G. - ESPINOSA, O. - MOLL, V.: {\em On some families of
integrals solvable in terms
of polygamma and negapolygamma functions}. To appear in {\em Integral 
Transforms and Special Functions}, 2003. 

\bibitem{bromwich}
BROMWICH, T.J.: {\em An Introduction to the Theory of Infinite Series}, 2nd.
Edition, MacMillan, New York, 1926.

\bibitem{esmo1}
ESPINOSA, O. - MOLL, V.: {\em On some integrals involving the Hurwitz
zeta function: part 1}. {The Ramanujan Journal}, {\bf 6}, 2002, 159-188.

\bibitem{esmo2}
ESPINOSA, O. - MOLL, V.: {\em On some integrals involving the Hurwitz
zeta function: part 2}. {The Ramanujan Journal}, {\bf 6}, 2002, 449-468.

\bibitem{gosper}
GOSPER, R. Wm. Jr.: $\int_{n/4}^{m/6} \ln \Gamma(z) dz$. In {\em Special
functions, $q$-series and related topics}, pages 71-76.
M. Ismail, D. Masson, M. Rahman
editors. The Fields Institute Communications, AMS, 1997.

\bibitem{gr}
GRADSHTEYN, I.S. - RYZHIK, I.M.: {\em Table of Integrals, Series and
Products}. Fifth edition, ed. Alan Jeffrey. Academic Press, 1994.

\bibitem{grossman}
GROSSMAN, N.: {\em Polygamma functions of arbitrary order}.
SIAM J. Math. Anal. {\bf 7}, 1976, 366-372.

\bibitem{knopp}
KNOPP, M. - ROBINS, S.: {\em Easy proofs of Riemann's functional equation
for} $\zeta(s)$ {\em and of Lipschitz summation}. Proc. AMS {\bf 129}, 2001,
1915-1922.

\bibitem{prudnikov}
PRUDNIKOV, A.P.- BRYCHKOV, Yu. A. - MARICHEV, O.I.: {\em Integrals and
Series}. Volume 3: {\em More special functions}. Translated from the
Russian by G.G. Gould. Gordon and Breach Science Publishers, New York, 1990.

\bibitem{ross}
ROSS, B.: {\em Problem 6002}.
Amer. Math. Monthly {\bf 81}, 1974, 1121.


\bibitem{ww}
WHITTAKER, E. - WATSON, G.: {\em A Course of Modern Analysis}. Cambridge
University Press, Fourth Edition reprinted, 1963.

\end{thebibliography}
\end{document}
</pre></TT></TD></TR></TABLE>